\documentclass[11pt]{article}

\usepackage[utf8]{inputenc}
\usepackage[T1]{fontenc}
\usepackage{mathpazo}

\RequirePackage[english]{babel}

\RequirePackage[a4paper,margin=2.5cm]{geometry}
\RequirePackage{parskip}

\newcommand{\affiliation}{\footnote}
\makeatletter
\def\@fnsymbol#1{\ensuremath{\ifcase#1\or *\or \dagger\or \ddagger\or \mathsection\or \|\or **\or \dagger\dagger \or \ddagger\ddagger \else\@ctrerr\fi}}
\makeatother

\usepackage{xcolor}
\definecolor{cblue}{RGB}{0,70,140}
\definecolor{cgreen}{RGB}{100,140,0}
\definecolor{cred}{RGB}{190,10,50}

\usepackage[shortlabels]{enumitem}
\setlist{itemsep=0ex,topsep=0ex,parsep=0.4ex}

\usepackage{amsmath,amsfonts,amssymb,amsthm,mathtools,bbm,centernot}

\usepackage[hyphens]{url} 
\usepackage[hidelinks,colorlinks,pagebackref]{hyperref} 
\hypersetup{colorlinks=true,citecolor=cyan,linkcolor=teal,urlcolor=cyan}
\usepackage[capitalise,nameinlink,noabbrev,compress]{cleveref} 

\renewcommand*{\backref}[1]{}
\renewcommand*{\backrefalt}[4]{
	\ifcase #1 Not cited.%
	\or $\uparrow$#2%
	\else $\uparrow$#2%
	\fi%
}

\let\OLDthebibliography\thebibliography
\renewcommand\thebibliography[1]{
	\OLDthebibliography{#1}
	\setlength{\parskip}{0pt}
	\setlength{\itemsep}{0pt plus 0.3ex}
}

\theoremstyle{plain}
\newtheorem{theorem}{Theorem}[section]

\newtheorem{claim}[theorem]{Claim}

\theoremstyle{definition}



\newcommand{\A}{{\cal A}}
\newcommand{\B}{{\cal B}}
\title{Largest density of a layered subgraph of a hypercube}
\author{Maria Axenovich\affiliation{Karlsruhe Institute of Technology, Karlsruhe, Germany, \href{mailto:maria.aksenovich@kit.edu}{\tt maria.aksenovich@kit.edu}.} \and Arsenii Sagdeev\affiliation{Karlsruhe Institute of Technology, Karlsruhe, Germany, \href{mailto:sagdeevarsenii@gmail.com}{\tt sagdeevarsenii@gmail.com}.}}

\date{\today}

\begin{document}
\maketitle

\begin{abstract}
Let $L(t)$ denote the largest number of edges induced by $t$ vertices from two vertex layers of a hypercube.  We show that
$$
\frac14 t\log_2 t+\frac18 t\log_2\log_2 t-O(t) \leq  L(t) \leq \frac14 t\log_2 t+ (1+o(1))t\log_2\log_2 t.$$
\end{abstract}

\section{Introduction}
For a natural number $n$, a hypercube $Q_n$ on a ground set $[n]=\{1, \ldots, n\}$ is a graph whose vertices are subsets of the ground set and two vertices are adjacent if and only if their symmetric difference has size $1$. Hypercubes are fundamental objects in geometry, combinatorics, graph theory, coding theory, and discrete optimization. They model highly symmetric communication networks, parallel-processing architectures, Boolean functions, and the structure of binary data.

A graph is called {\it cubical} if it is isomorphic to a subgraph of a hypercube. An {\it edge-layer} of $Q_n$ is a subgraph induced by $\binom{[n]}{k} \cup \binom{[n]}{k-1}$, for some $k\in [n]$.
A graph is called {\it layered} if it is isomorphic to a subgraph of some edge-layer. For example, a cycle on $6$ vertices is layered, but a cycle on $4$ vertices is not. All logarithms in this paper are base 2.

\begin{theorem}\label{main}
Let $L(t)$ be the largest number of edges in a $t$-vertex layered graph.  Then 
$$
\frac14 t\log t+\frac18 t\log\log t-O(t) \leq  L(t) \leq \frac14 t\log t+ (1+o(1))t\log\log t.$$
\end{theorem}

Note that the maximal edge-density of a layered graph is at most the maximal density of a general cubical graph. The corresponding function $C(t)$, equal to the largest number of edges in a $t$-vertex cubical graph, has been studied in the context of isoperimetric problems, see for example Leader~\cite{L}.
Harper~\cite{H}  determined that 
$C(t)= \sum_{i=0}^{t-1} s_i$, where $s_i$ is the number of $1$'s in the binary representation of $i$. See also Graham~\cite{G},  Bollob\'as~\cite{B}, and Hart \cite{Ha}. The Trollope–Delange formula gives $C(t)=\frac{1}{2}t\log t + t\mathcal{T}(\log t)$ for some continuous periodic function $\mathcal{T}$, see Trollope~\cite{T}, Delange~\cite{D}, and also Flajolet, Grabner, Kirschenhofer, Prodinger and Tichy~\cite{F}. In particular, there is a positive constant $c$ such that
$$\frac{1}{2}t\log t -ct\leq C(t) \leq \frac{1}{2}t\log t +ct.$$
Thus, Theorem~\ref{main} implies that $L(t)=(\frac12+o(1)) C(t)$.

 For a graph $H$, let ${\rm ex}(Q_n, H)$ be the largest number of edges in a subgraph of $Q_n$ without a subgraph isomorphic to $H$. Identifying graphs having positive Tur\'an density in the hypercube, i.e., those $H$ for which ${\rm ex}(Q_n, H)= \Omega(|E(Q_n)|)$ remains an open problem. See  for example Conlon \cite{C}, Axenovich~\cite{A}, Axenovich, Martin, Winter~\cite{AMW}, and  Zhu~\cite{Zh}.  
 One can see that for every non-layered graph $H$, ${\rm ex}(Q_n, H) \geq \frac{1}{2} |E(Q_n)|$ by taking a union of every other edge-layer of $Q_n$. However, there were no quantitative condition to identify cubical non-layered graphs.
 Theorem~\ref{main} sheds more light on this problem, since any graph $H$ on $t$ vertices and more than $L(t)$ edges is not layered.



\section{The upper bound}

Let $L(t)$ be the largest number of edges induced by $t$ vertices from two vertex layers of a hypercube.  Define $\Phi(0)=0$ and, for $t\ge1$,
$$ \Phi(t)=\frac14t\log t+t\log(\log t+2).$$
We shall show by induction on $t$ that 
$L(t)\leq \Phi(t).$
The base cases $L(0)=0=\Phi(0)$ and $L(1)=0<1=\Phi(1)$ are immediate. In what follows, we assume that $t \ge 2$.

Consider a set $\cal A \cup \cal B$ of $t$ vertices in a hypercube  on the ground set $[n]$, where each vertex in $\cal A$ has size $k$ and each vertex in $\cal B$ has size $k+1$, for some $k$ and $n$.  Assume that $\A\cup \B$ induces the largest possible number $L(t)$ of edges. In addition, assume that $n$ is minimal among all such examples.
For sets $\A$ and $\cal B$ of vertices in a given graph, let $e(\A, \B)$ be the number of edges with one endpoint in $\cal A$ and the other in $\cal B$.  Next, we shall upper bound $L(t)=e(\A, \B)$.

For $i\in [n]$, let $\A\cup \B= \A_0 (i)\cup \A_1(i)\cup \B_0(i)\cup \B_1(i),$ where
$$\A_0 (i)=\{X\in \A: i\not\in X\}, ~~\A_1(i)= \{X\in \A: i\in X\},$$
$$\B_0(i)= \{X\in \B: i\not\in X\}, ~~\B_1(i)= \{X\in \B: i\in X\}.$$
Since there are no edges between $\B_0(i)$ and $\A_1(i)$,  and the edges between 
$\B_1(i)$ and $\A _0(i)$ form a matching, we have
 that $$e(\A, \B)\leq  e(\A_0(i), \B_0(i))+ e(\A_1(i), \B_1(i)) + \min \{|\A_0(i)|, |\B_1(i)|\}.$$
Observe that $|\A_0 (i)\cup \B_0(i)| < |\A\cup \B|$ by the minimality of $n$ since otherwise we could remove $i$ from the ground set to get a smaller example. Thus $e(\A_0(i), \B_0(i)) \le L(|\A_0 (i)\cup \B_0(i)|) \le \Phi(|\A_0 (i)\cup \B_0(i)|)$ by induction. Similarly, $|\A_1 (i)\cup \B_1(i)| < |\A\cup \B|$ and thus $e(\A_1(i), \B_1(i)) \le \Phi(|\A_1 (i)\cup \B_1(i)|)$. Using the notation  $a_0(i):= |\A_0(i)|/t$,~ $a_1(i):= |\A_1(i)|/t$,~ $b_0(i):= |\B_0(i)|/t$, $b_1(i):= |\B_1(i)|/t$, and $r_i=a_0(i)+ b_0(i)$, we get
 \begin{equation}\label{eq-main}
 e(\A, \B)\leq  \Phi( r_i t) + \Phi((1-r_i)t) + t \min \{a_0(i), b_1(i)\}.
 \end{equation}

We shall find a ``good'' element $i$ such that the term $\min \{a_0(i), b_1(i)\}$ can be upper bounded by a convenient function. For that we shall use the binary entropy function $H(x)$ defined for $0\le x\le1$, as follows
 $H(x)=-x\log x-(1-x)\log(1-x),$ with the convention $0\log0=0$. 

\begin{claim} \label{cl1}
    There is $i\in [n]$ such that 
$ \min\{a_0(i), b_1(i)\} \leq \frac{1}{4}H(r_i) + {H(r_i)}/{(4\log t)}.$
\end{claim}

To prove the claim, let $a=|\A|/t= a_0(i)+a_1(i)$. Recall that $1= a_0(i)+a_1(i)+b_0(i)+b_1(i)$ and $r_i= a_0(i)+b_0(i)$. 
We define $\Delta_i = a_0(i) - ar_i.$
Then $b_1(i)= (1-a)(1-r_i)+\Delta_i$
and
$$ \min\{a_0(i), b_1(i)\} = \min\{ar_i,(1-a)(1-r_i)\} +\Delta_i\le r_i(1-r_i) +\Delta_i \le \frac14H(r_i)+ \Delta_i. $$
In the last inequality we used the standard fact that $H(x)\ge 4x(1-x)$, see Claim~\ref{cl5} in Appendix.
Thus, in order to prove the claim it is enough to find an $i$ for which $\Delta_i\le H(r_i)/(4\log t)$.

Since the members of $\A$ are sets of size $k$ and the members of $\B$ are sets of size $k+1$ each, we have 
$
 \sum_{i\in [n]} a_1(i)=ka,~
 \sum_{i\in [n]} b_1(i)=(k+1)(1-a).$
Since 
$ \Delta_i=ab_1(i)-(1-a)a_1(i),$
we have 
$$
 \sum_{i\in [n]}\Delta_i=a(1-a)\le\frac14.
$$
By implicit usage of entropy subadditivity, see Claim~\ref{cl6} in Appendix, we have
$$ \sum_{i\in [n]}\frac{H(r_i)}{4\log t}\ge \frac14.
$$

Thus there is an $i\in [n]$ such that $\Delta_i\leq H(r_i)/(4\log t)$. 
This proves Claim~\ref{cl1}.\\

Recall that  $\Phi(0)=0$ and, for $t\ge1$, $ \Phi(t)=\frac14t\log t+t\log(\log t+2).$

\begin{claim} \label{cl2}
    Let  $t\ge2$ and $0< r< 1$, such that $rt$ and $(1-r)t$ are integers. Then 
 $$\Phi(rt)+\Phi((1-r)t)
 \le
\Phi(t) - \frac{t}{4}H(r)-\frac{tH(r)}{4\log t}.$$
\end{claim}

To prove the claim, observe first that
$\frac14t\log t- \frac14 rt\log (rt) -\frac14(1-r)t\log ((1-r)t) = \frac{t}{4}H(r). $ It remains to estimate the terms involving $G(t)=t\log(\log t+2)$:
\begin{align*}
G(t)-G(rt)-G((1-r)t)
&=rt\log\frac{\log t+2}{\log(rt)+2}
 +(1-r)t\log\frac{\log t+2}{\log((1-r)t)+2}  \\
&=-rt\log\left(1+\frac{\log r}{\log t+2}\right)
 -(1-r)t\log\left(1+\frac{\log (1-r)}{\log t+2}\right)  \\
&\ge \frac{t}{(\log t+2)\ln2}
 \left(-r\log r-(1-r)\log(1-r)\right), \\
&\ge \frac{tH(r)}{4\log t},
\end{align*}
as desired. Here we used the inequalities $-\log(1+y)\ge -y/\ln2$ for $-1<y<0$ and $1/{((\log t+2)\ln2)}\ge 1/{(4\log t)}$ for $t\ge2$.
This proves Claim~\ref{cl2}.\\

Now, we return to inequality (\ref{eq-main}) and use Claims~\ref{cl1} and~\ref{cl2}:
\begin{eqnarray*}
 L(t)=e(\A, \B) &\leq & \Phi( r_i t) + \Phi((1-r_i)t) + t \min \{a_0(i), b_1(i)\}\\
   &\leq & \Phi( r_i t) + \Phi((1-r_i)t)  + \frac{t}{4}H(r_i) + \frac{t H(r_i)}{4\log t}\\
&\leq & \Phi(t).   
\end{eqnarray*}
This proves the upper bound of Theorem~\ref{main}.

\section{The lower bound}

For $m \ge 1$, let  \(v_m\) and \(e_m\) be the numbers of vertices and edges,
respectively, in the middle edge-layer of \(Q_m\), i.e., for 
$
 k=\left\lfloor\frac{m-1}{2}\right\rfloor,$ 
 $v_m=\binom mk+\binom m{k+1}
     =\binom{m+1}{k+1}$ and $
 e_m=(k+1)\binom m{k+1}.$

\begin{claim} \label{cl3}
    We have $e_m=\frac14v_m\log v_m
     +\frac18v_m\log\log v_m+O(v_m).$
Moreover, \(v_m <v_{m+1}\le2v_m\).
\end{claim}

Since \(v_m\) is a middle binomial coefficient, Stirling's
formula gives
$ \log v_m=m-\frac12\log m+O(1).$
Also,
$ \frac{e_m}{v_m}=\frac m4+O(1).$
Hence $m=\log v_m+\frac12\log\log v_m+O(1),$
and the asserted estimate for \(e_m\) follows.
For the second assertion, Pascal's identity gives
$ v_{2\ell}=\binom{2\ell+1}{\ell},$ $v_{2\ell+1}=\binom{2\ell+2}{\ell+1}=2v_{2\ell}$,
and
$v_{2\ell+2}=\binom{2\ell+3}{\ell+1}
 <2\binom{2\ell+2}{\ell+1}=2v_{2\ell+1}.$
Thus \(v_m < v_{m+1}\le2v_m\).
This proves Claim~\ref{cl3}.\\

We now define \(F(0)=0\) and, for \(t\ge1\),
\[
 F(t)=\frac14t\log t+\frac18t\log(\log t+2).
\]
By Claim~\ref{cl3}, \(e_m\ge F(v_m)-C_0v_m\), for some absolute constant \(C_0\ge0\). Let $C=C_0+\frac{3}{4}$.
We shall show by induction on $t$ that \(L(t)\ge F(t)-Ct\). The base cases $t=0,1$ are immediate. In what follows, we assume that $t \ge 2$.

Choose \(m\) maximal such that \(v_m\le t\), and let \(r=v_m/t\). By the maximality of \(m\) and Claim~\ref{cl3}, \(t<v_{m+1}\le2v_m=2rt\), so \(r>\frac{1}{2}\).  If $r=1$, then we are done since \(C>C_0\), so assume that $r<1$. The induction hypothesis gives a layered graph on \((1-r)t\) vertices with at least \(F((1-r)t)-C(1-r)t\) edges. Taking its vertex-disjoint union with the middle edge-layer of \(Q_m\), we obtain a \(t\)-vertex layered graph. Hence
\begin{equation}\label{eq-lb}
 L(t)\ge F(rt)-C_0rt+F((1-r)t)-C(1-r)t .
\end{equation}

\begin{claim} \label{cl4}
Let  $t\ge2$ and $0< r< 1$, such that $rt$ and $(1-r)t$ are integers. Then 
 $$F(rt)+F((1-r)t) \ge
F(t) - \frac{3}{8}t.$$
\end{claim}

To prove the claim, observe first that
\(\frac14t\log t-\frac14rt\log(rt)-\frac14(1-r)t\log((1-r)t)=\frac{t}{4}H(r)\le \frac{t}{4}.\)
It remains to estimate the terms involving \(G(t)=t\log(\log t+2)\):
\[
\begin{aligned}
G(t)-G(rt)-G((1-r)t)
&=rt\log\frac{\log t+2}{\log(rt)+2}
 +(1-r)t\log\frac{\log t+2}{\log((1-r)t)+2}  \\
&\le rt\log\frac{t}{tr}+(1-r)t\log\frac{t}{(1-r)t} \\
&=tH(r) \le t,
\end{aligned}
\]
where we used that the function \((\log t+2)/t\) is decreasing for $t \ge 1$. This proves Claim~\ref{cl4}. \\

Recall that $r > \frac{1}{2}$ and $C=C_0+\frac{3}{4}$. Now, we return to inequality (\ref{eq-lb}) and use Claim~\ref{cl4}:
\[L(t) \ge  F(t)-\frac38t-C_0rt-C(1-r)t = F(t)-Ct +\left(\frac{3}{4}r-\frac38\right)t \ge F(t)-Ct.\]
This establishes the lower bound in Theorem~\ref{main}. \\


{\bf Remark}\quad A more careful analysis of Claim~\ref{cl2} shows that for any $c>\frac{\ln2}{4}$, we have $L(t)\le \frac14 t\log t+ ct\log\log t+O_c(t)$.  ChatGPT-5.5 Pro provided a more involved upper bound proof that gives $c=\frac{1}{8}$ matching the lower bound up to an $o(t \log \log t)$ term. It would be interesting to find a closed-form expression for $L(t)$, similar to Harper's formula for $C(t)$.


{\bf Acknowledgements}  The proof of the main theorem was originally generated with the assistance of ChatGPT-5.5 Pro, then checked, shortened, and reworked by the authors.

\appendix

\section{Properties of the entropy function} \label{entropy}

\begin{claim} \label{cl5}
    For all $0\le x \le 1$, we have $f(x) = H(x)-4x(1-x) \ge 0$.
\end{claim}

Since $f(x)=f(1-x)$, we can assume that \(0\le x\le \frac12\). We have
\(
f''(x)=8-\frac{1}{(\ln 2)x(1-x)},
\)
so \(f''\) is strictly increasing and changes sign exactly once on \(\left(0,\frac12\right)\). 
Moreover, \(
f(0)=f\left(\frac12\right)=0,
\)
\( f'\left(\frac12\right)=0,\) and $f'(x)>0$ for sufficiently small $x>0$.
Hence \(f'\) first decreases and then increases, so it changes sign exactly once on \(\left(0,\frac12\right)\). Therefore \(f\) first increases and then decreases back to \(0\), and thus
\(
f(x)\ge 0
\)
for \(0\le x\le \frac12\), as claimed.

\begin{claim} \label{cl6}
    $\sum_{i=1}^n H(r_i)\ge \log t.$
\end{claim}

Recall that 
$\mathcal F=\mathcal A\cup\mathcal B,$ $|\mathcal F|=t,$
and $r_i
={|\{X\in\mathcal F:i\notin X\}|}/{t}.$
For each $X\subseteq[n]$, define
$
w(X)=\big(\prod_{i\notin X}r_i\big)\big(
     \prod_{i\in X}(1-r_i)\big).$
Then
$
\sum_{X\subseteq[n]}w(X)
=
\prod_{i=1}^n\bigl(r_i+(1-r_i)\bigr)
=1.$
In particular,
$
\sum_{X\in\mathcal F}w(X)\le 1.$
Hence $ \left(\prod_{X\in\mathcal F}w(X)\right)^{1/t} \le (\sum_{X\in\mathcal F}w(X))/t \le 1/t$ by the arithmetic--geometric mean inequality.
Now coordinate \(i\) is absent from exactly \(r_it\) members of
\(\mathcal F\), and belongs to exactly \((1-r_i)t\) members.
Therefore
$
\prod_{X\in\mathcal F}w(X)
=
\prod_{i=1}^n
r_i^{r_it}(1-r_i)^{(1-r_i)t}.$
Taking \(t\)-th roots, we obtain
$\prod_{i=1}^n
r_i^{r_i}(1-r_i)^{1-r_i}
\le 1/t.$
Taking logarithms gives
$\sum_{i=1}^n H(r_i)
\ge \log t,$ as claimed.

\end{document}